\documentclass[a4paper,12pt]{article}
\usepackage{stmaryrd}
\usepackage{amssymb,amsthm,amsmath,amsfonts}
\usepackage[dvips]{graphicx}
\usepackage{geometry}

\begin{document}

\newtheorem{thm}{Theorem}[]
\newtheorem{pro}{Proposition}[]
\newtheorem{lem}{Lemma}[]

\title{The Second Variation of the Ricci Expander Entropy}
\author{Meng Zhu\thanks{Research is partially supported by NSF grant DMS-0354621.}}
\date{}

\maketitle

\section{Introduction}

\hspace{0.5cm} In \cite{Pe}, Perelman discovered two important
functionals, the $\mathcal{F}$-functional and the $\mathcal{W}$
functional. The corresponding entropy functionals $\lambda$ and
$\nu$ are monotone along the Ricci flow $\frac{\partial
g_{ij}}{\partial t}=-2R_{ij}$ and constant precisely on steady and
shrinking solitons. In \cite{Cao}, H.-D. Cao, R. Hamilton and T.
Ilmanen presented the second variations of both entropy functionals and
studied the linear stabilities of certain closed
Einstein manifolds of nonnegative scalar curvature.\\

To find the corresponding variational structure for the expanding
case, M. Feldman, T. Ilmanen and L. Ni \cite {Ni} introduced the
$\mathcal{W}_+$ functional. Let $(M^n,g)$ be a compact Riemannian
manifold, $f$ a smooth function on $M$, and $\sigma>0$. Define
$$\mathcal{W}_+(g,f,\sigma)=(4\pi \sigma)^{-\frac{n}{2}}\int_M e^{-f}[\sigma(|\nabla f|^2 + R)-f+n]dV,$$
$$\mu_+(g,\sigma)=\inf\{W_+(g,f,\sigma)| \quad f\in C^{\infty}(M), \quad and \quad (4\pi \sigma)^{-\frac{n}{2}}\int_M e^{-f}dV=1\},$$
and
$$\nu_+(g)=\sup_{\sigma >\ 0}\mu_+(g,\sigma).$$
Then $\nu_+$ is nondecreasing along
the Ricci flow and constant precisely on expanding solitons.\\

In this note, analogous to \cite{Cao}, we present the first and
second variations of the entropy $\nu_+$. By computing the first
variation of $\nu_+$, one can see that the critical points are
expanding solitons, which are actually
negative Einstein manifolds (see e.g. \cite{Cao2}). Our main result is the following\\

\begin{thm}
Let $(M^n,g)$ be a compact negative Einstein manifold. Let $h$ be a
symmetric 2-tensor. Consider the variation of metric $g(s)=g+sh$.
Then the second variation of $\nu_+$ is

$$\frac{\mathrm{d}^2 \nu_+(g(s))}{\mathrm{d}s^2}|_{s=0} =
\frac{\sigma}{Vol(g)}\int_M <N_+h,h>,$$

\hspace{-0.6cm}where $$N_+h:=\frac{1}{2}\Delta h +
\mathrm{div^*div}h + \frac{1}{2}\nabla^2 v_h + Rm(h,\cdot)+
\frac{g}{2n\sigma \mathrm{vol}(g)}\int_M trh,$$

\hspace{-0.6cm}and $v_h$ is the unique solution of
$$\Delta v_h - \frac{v_h}{2\sigma}= \mathrm{div(div}h), \hspace{1cm} \int_M v_h=0.$$
\end{thm}

In this case, we may still define the concept of linear stability.
We say that an expanding soliton is \textit{linearly stable} if
$N_+\leq 0$, otherwise it is \textit{linearly unstable}. Similar to
\cite{Cao}, the $N_+$ operator is nonpositive definite if and only
if the maximal eigenvalue of the Lichnerowicz Laplacian acting on
the space of transverse traceless 2-tensors has certain upper bound.
Then using the results in \cite{De} and \cite{De2}, one can see that
compact hyperbolic spaces are linearly stable. But unlike the
positive Einstein case, it seems hard to find other examples of
negative Einstein manifolds which are either linear stable or
linear unstable. \\

\hspace{-.75cm} \textbf{Acknowledgement} The author would like to
thank Professor Huai-Dong Cao for his encouragement and
suggestions.\\

\section{\textbf{The First Variation of the Expander Entropy}}

Recall that in \cite{Pe}, the $\mathcal{F}$ functional is defined by

$$\mathcal{F}(f,g)=\int_M(|\nabla f|^2 + R)e^{-f}\mathrm{d}V,$$

\hspace{-0.6cm}and its entropy $\lambda(g)$ is

$$\lambda(g)=\inf\{\mathcal{F}(f,g): f\in C^{\infty}(M), \int_M e^{-f}=1\},$$

\hspace{-0.6cm}where $R$ is the scalar curvature. By Theorem 1.7 in
\cite{Ni}, we know that $\mu_+(g,\sigma)$ is attained by some
function $f$. Moreover, if $\lambda(g)<0$, then $\nu_+(g)$ can be
attained by some positive number $\sigma$.

\begin{lem}
Assume that $\nu_+(g)$ is realized by some $f$ and $\sigma$, then it
is necessary that the pair $(f,\sigma)$ solves the following
equations,

\begin{equation}
\sigma(-2\Delta f + |\nabla f|^2 - R)+f-n+\nu_+=0,
\end{equation}
and
\begin{equation}
(4\pi\sigma)^{-\frac{n}{2}}\int_M fe^{-f}dV=\frac{n}{2}-\nu_+.
\end{equation}
\end{lem}

\textbf{Proof}: For fixed $\sigma>0$, suppose that $\mu_+(g,\sigma)$
is attained by some function $f$. Using Lagrange multiplier method,
consider the following functional

$$L(g,f,\sigma,\lambda)=(4\pi \sigma)^{-\frac{n}{2}}\int_M e^{-f}[\sigma(|\nabla f|^2 + R)-f+n]dV-\lambda((4\pi \sigma)^{-\frac{n}{2}}\int_M e^{-f}dV-1).$$

Denote by $\delta f$ the variation of $f$. Then the variation of $L$
is

\begin{eqnarray*}
0=\delta L &=& (4\pi \sigma)^{-\frac{n}{2}}\int_M e^{-f}\{(-\delta
f)[\sigma(|\nabla f|^2 + R)-f+n]+[2\nabla f \nabla (\delta f)-\delta
f]\}\mathrm{d}V\\
& & \ - (4\pi \sigma)^{-\frac{n}{2}}\int_M \lambda (\delta
f) e^{-f}\mathrm{d}V\\
& = & (4\pi \sigma)^{-\frac{n}{2}}\int_M e^{-f}(\delta
f)[\sigma(-2\Delta f+|\nabla f|^2-R)+f-n-1-\lambda]\mathrm{d}V
\end{eqnarray*}

Therefore, we have

$$\sigma(-2\Delta f+|\nabla f|^2-R)+f-n-1-\lambda=0.$$

Integrating both sides with respect to the measure $(4\pi
\sigma)^{-\frac{n}{2}}e^{-f}\mathrm{d}V$, we get

$$-\lambda - 1= (4\pi \sigma)^{-\frac{n}{2}}\int_M e^{-f}[\sigma(|\nabla f|^2 + R)-f+n]dV= \mu_+(g,\sigma).$$

When $\sigma$ and $f$ realize $\nu_+(g)$, the above formula is just
equation $(1)$.\\

Now we consider the variations $\delta \sigma$ and $\delta f$ of
both $\sigma$ and $f$. We have

\parbox{10cm}
{\begin{eqnarray*} 0 &=& (4\pi \sigma)^{-\frac{n}{2}}\int_M
e^{-f}(-\frac{n}{2\sigma}\delta \sigma-\delta f)[\sigma(|\nabla f|^2
+ R)-f+n]\mathrm{d}V\\
& & \ +\ (4\pi \sigma)^{-\frac{n}{2}}\int_M e^{-f}[\delta
\sigma(|\nabla f|^2+R)+ 2\nabla f \nabla (\delta f)-\delta
f]\mathrm{d}V
\end{eqnarray*}}\hfill
\parbox{1cm}{\begin{eqnarray}\end{eqnarray}}

and

\begin{equation}
(4\pi \sigma)^{-\frac{n}{2}}\int_M e^{-f}(-\frac{n}{2\sigma}\delta
\sigma-\delta f)\mathrm{d}V=0.
\end{equation}

Using $(1)$ and $(4)$, we can write $(3)$ as

\begin{eqnarray*}
0 & =& (4\pi \sigma)^{-\frac{n}{2}}\int_M
e^{-f}[\delta\sigma(|\nabla
f|^2 + R)-\delta f]\mathrm{d}V\\
&=& (4\pi \sigma)^{-\frac{n}{2}}\int_M
e^{-f}[\frac{1}{\sigma}\delta\sigma(\nu_+
+f-n)+\frac{n}{2\sigma}\delta\sigma]\mathrm{d}V\\
& = & (\delta \sigma) \frac{1}{\sigma}(4\pi
\sigma)^{-\frac{n}{2}}\int_M e^{-f}(\nu_+ +f-\frac{n}{2})\mathrm{d}V
\end{eqnarray*}

Hence, we get equation (2). \hspace{7.5cm} Q.E.D.\\

Before computing the variations of $\nu_+$ functional, let's recall
some variation formulas of curvatures.  By direct computation, we
can get the following lemma

\begin{lem}
Suppose that $h$ is a symmetric 2-tensor, and $g(s)=g+sh$ is a
variation of $g$. Then

\begin{equation}
\frac{\partial R}{\partial s}|_{s=0}= -h_{kl}R_{kl}+\nabla_p\nabla_k
h_{pk}- \Delta trh,
\end{equation}

and

\parbox{10cm}
{\begin{eqnarray*} \frac{\partial^2 R}{\partial s^2}|_{s=0} & = &
2h_{kp}h_{pl}R_{kl} - 2h_{kl}\frac{\partial R_{kl}}{\partial
s}|_{s=0} +
g^{kl}\frac{\partial^2 R_{kl}}{\partial s^2}|_{s=0} \\
& = & 2h_{kp}h_{pl}R_{kl} - h_{kl}(2\nabla_p\nabla_k
h_{pl}-\Delta h_{kl}- \nabla_k\nabla_l \mathrm{tr}h)\\
& &  -\  \nabla_p [h_{pq}(2\nabla_k h_{kq}- \nabla_q \mathrm{tr}h)] + \nabla_k (h_{pq}\nabla_k h_{pq})\\
& & +\  \frac{1}{2}\nabla_p \mathrm{tr}h(2\nabla_k h_{kp} -\nabla_p \mathrm{tr}h)\\
& & +\   \frac{1}{2}(\nabla_k h_{pq} \nabla_k h_{pq}- 2\nabla_p
h_{kq}\nabla_q h_{kp})
\end{eqnarray*}}\hfill
\parbox{2.5cm}{\begin{eqnarray}\end{eqnarray}}

Here, $\nabla$ is the Levi-Civita connection of $g$, and
$\mathrm{tr}h$ is the trace of $h$ taken with respect to $g$.\\
\end{lem}

Now we are ready to compute the first variation of
$\nu_+(g)$.\\

\begin{pro}
Let $(M^n,g)$ be a compact Riemannian manifold with $\lambda(g)<0$.
Let $h$ be any symmetric covariant 2-tensor on $M$, and consider the
variation $g(s)=g+sh$. Then the first variation of $\nu_+(g(s))$ is

$$\frac{d\nu_+(g(s))}{ds}|_{s=0}= (4\pi \sigma)^{-\frac{n}{2}}\int_M \sigma e^{-f}(-R_{ij}-\nabla_i\nabla_jf -\frac{1}{2\sigma}g_{ij})h_{ij}dV,$$
where the smooth function $f$ and $\sigma>0$ realize $\nu_+(g)$. \\
\end{pro}

\textbf{Proof}:

\parbox{10cm}
{\begin{eqnarray*} \frac{\partial \nu_+}{\partial s} & = &
(4\pi\sigma)^{-\frac{n}{2}} \int_M
e^{-f}(-\frac{n}{2\sigma}\frac{\partial \sigma}{\partial s}
-\frac{\partial f}{\partial s} +
\frac{1}{2}g^{ij}h_{ij})[\sigma(|\nabla f|^2+R)-f+n]dV\\
& & +\  (4\pi\sigma)^{-\frac{n}{2}}\int_M e^{-f}\frac{\partial
\sigma}{\partial s}(|\nabla f|^2+R)\mathrm{d}V\\
& & +\ (4\pi\sigma)^{-\frac{n}{2}}\int_M
e^{-f}[\sigma(-g^{ip}g^{jq}h_{pq}\nabla_if\nabla_jf+\
2g^{ij}\nabla_if\nabla_j\frac{\partial f}{\partial s}+\frac{\partial
R}{\partial s})-
\frac{\partial f}{\partial s}]dV.\\
\end{eqnarray*}}\hfill
\parbox{1cm}{\begin{eqnarray}\end{eqnarray}}

From
$$(4\pi\sigma)^{-\frac{n}{2}}\int_M e^{-f}\mathrm{d}V=1,$$
we have
\begin{equation}
(4\pi\sigma)^{-\frac{n}{2}}\int_M (-\frac{n}{2\sigma}\frac{\partial
\sigma}{\partial s}-\frac{\partial f}{\partial
s}+\frac{1}{2}g^{ij}h_{ij})e^{-f}\mathrm{d}V=0.
\end{equation}\\

Substituting $(1),\ (2)$ and $(8)$ in $(7)$, we obtain
\begin{eqnarray*} \frac{\partial \nu_+(s)}{\partial s}|_{s=0} & = &
(4\pi\sigma)^{-\frac{n}{2}}\int_M [2\sigma(|\nabla f|^2-\Delta
f)+\nu_+(0)](-\frac{n}{2\sigma}\frac{\partial \sigma}{\partial
s}-\frac{\partial f}{\partial
s}+\frac{1}{2}g^{ij}h_{ij})e^{-f}\mathrm{d}V\\
& & +\ (4\pi\sigma)^{-\frac{n}{2}}\int_M[\frac{\partial
\sigma}{\partial s}(|\nabla f|^2+R)-\frac{\partial f}{\partial s}+
\sigma(-h_{ij}\nabla_i f \nabla_j f + 2\frac{\partial f}{\partial
s}(|\nabla f|^2-\Delta f)\\
& & \hspace{3cm} +\  \nabla_i \nabla _j h_{ij}-\Delta trh -
h_{ij}R_{ij})]e^{-f}\mathrm{d}V\\
\hspace{2cm}& = & (4\pi\sigma)^{-\frac{n}{2}}\int_M [\frac{\partial
\sigma}{\partial s}(|\nabla f|^2+R)-\frac{\partial f}{\partial
s}-\sigma( h_{ij}\nabla_i\nabla_j f + h_{ij}R_{ij})]e^{-f}\mathrm{d}V\\
& = & (4\pi\sigma)^{-\frac{n}{2}}\int_M [\frac{\partial
\sigma}{\partial s}(|\nabla f|^2+R)+ \frac{n}{2\sigma}\frac{\partial
\sigma}{\partial s}-\sigma h_{ij}(R_{ij}+ \nabla_i\nabla_j f + \frac{1}{2\sigma}g_{ij})]e^{-f}\mathrm{d}V\\
& = & (4\pi\sigma)^{-\frac{n}{2}}\int_M
\frac{1}{\sigma}\frac{\partial \sigma}{\partial
s}[f(0)-\frac{n}{2}+\nu_+(0)-2\sigma(|\nabla f|^2-\Delta f)]e^{-f}\\
& &\hspace{2cm} -\ \sigma h_{ij}(R_{ij}+ \nabla_i\nabla_j f +
\frac{1}{2\sigma}g_{ij})e^{-f}\mathrm{d}V\\
& = & (4\pi\sigma)^{-\frac{n}{2}}\int_M -\sigma h_{ij}(R_{ij}+
\nabla_i\nabla_j f + \frac{1}{2\sigma}g_{ij})e^{-f}\mathrm{d}V.
\end{eqnarray*}

Hence, the first variation of $\nu_+$ is

$$\frac{d\nu_+(g(s))}{ds}|_{s=0}= (4\pi \sigma)^{-\frac{n}{2}}\int_M \sigma e^{-f}(-R_{ij}-\nabla_i\nabla_jf -
\frac{1}{2\sigma}g_{ij})h_{ij}dV.$$

\hspace{12.5cm}Q.E.D.\\

From the above proposition, we can see that a critical point of
$\nu_+(g)$ satisfies

$$Rc+\nabla^2f+\frac{1}{2\sigma}g=0,$$

\hspace{-0.6cm}which means that $(M,g)$ is a gradient
expanding soliton.\\

\section{The Second Variation}

\hspace{.6cm}Now we compute the second variation of $\nu_+$. Since
any compact expanding soliton is Einstein (e.g. see \cite{Cao2}), it
implies that $f$ is a constant. After adding $f$ by
a constant we may assume that $f=\frac{n}{2}$.\\

In the following, as in \cite{Cao}, we denote
$Rm(h,h)=R_{ijkl}h_{ik}h_{jl}$,
$\mathrm{div}\omega=\nabla_i\omega_i$,
$(\mathrm{div}h)_i=\nabla_jh_{ji}$,
$(\mathrm{div^*}\omega)_{ij}=-(\nabla_i\omega_j+\nabla_j\omega_i)=-\frac{1}{2}L_{\omega^\#}g_{ij}$,
where $h$ is a symmetric 2-tensor, $\omega$ is a 1-tensor,
$\omega^{\#}$ is the dual vector field of $\omega$, and
$L_{\omega^{\#}}$ is the Lie derivative.\\

\textbf{Proof of Theorem 1}:  Let $(M,g)$ be a compact negative
Einstein manifold with $f=\frac{n}{2}$ and
$R_{ij}=-\frac{1}{2\sigma}g_{ij}$. For any symmetric 2-tensor $h$,
consider the variation $g(s)=g+sh$. Then by proposition 1, we know
that $\frac{d\nu_+}{ds}|_{s=0}=0$.\\

From (1) and (2), we can get

\begin{equation}
\frac{n}{2\sigma}\frac{\partial \sigma}{\partial
s}(0)-2\sigma\Delta\frac{\partial f}{\partial s}(0)-\sigma
\frac{\partial R}{\partial s}(0)+ \frac{\partial f}{\partial
s}(0)=0,
\end{equation}

\hspace{-0.6cm}and

$$(4\pi\sigma)^{-\frac{n}{2}}\int_M \frac{n}{2}e^{-\frac{n}{2}}(-\frac{n}{2\sigma}\frac{\partial \sigma}{\partial s}(0)-\frac{\partial f}{\partial s}(0)+ \frac{1}{2}tr_gh)dV + (4\pi\sigma)^{-\frac{n}{2}}\int_M \frac{\partial f}{\partial s}(0)e^{-\frac{n}{2}}dV=0.$$\\

It follows by (8) that

\begin{equation}
(4\pi\sigma)^{-\frac{n}{2}}\int_M \frac{\partial f}{\partial
s}(0)e^{-\frac{n}{2}}dV=0,
\end{equation}
and
\begin{equation}
\frac{n}{2\sigma}\frac{\partial \sigma}{\partial s}(0)=
\frac{1}{Vol(g)}\int_M \frac{1}{2}trhdV,
\end{equation}

\hspace{-0.6cm}where
$(4\pi\sigma)^{-\frac{n}{2}}e^{-\frac{n}{2}}=\frac{1}{Vol(g)}.$ Thus

\begin{eqnarray*}
\frac{\mathrm{d} \nu_+}{\mathrm{d} s} & = &
(4\pi\sigma)^{-\frac{n}{2}} \int_M
e^{-f}(-\frac{n}{2\sigma}\frac{\partial \sigma}{\partial s}
-\frac{\partial f}{\partial s} +
\frac{1}{2}g^{ij}h_{ij})[\sigma(|\nabla f|^2+R)-f+n]dV\\
& & +\  (4\pi\sigma)^{-\frac{n}{2}}\int_M e^{-f}[\frac{\partial
\sigma}{\partial s}(|\nabla f|^2+R)+
\sigma(-g^{ip}g^{jq}h_{pq}\nabla_if\nabla_jf+
2g^{ij}\nabla_if\nabla_j\frac{\partial f}{\partial s}+\frac{\partial
R}{\partial s})\\
& &\quad \quad \quad \quad \quad \quad \quad \quad  -\
\frac{\partial f}{\partial s}]dV\\
& = & (4\pi\sigma)^{-\frac{n}{2}} \int_M
e^{-f}(-\frac{n}{2\sigma}\frac{\partial \sigma}{\partial s}
-\frac{\partial f}{\partial s} +
\frac{1}{2}g^{ij}h_{ij})[2\sigma(|\nabla f|^2-\Delta f)+ \nu_+]dV\\
& & + \ (4\pi\sigma)^{-\frac{n}{2}}\int_M e^{-f}[\frac{\partial
\sigma}{\partial s}(|\nabla f|^2+R)+
\sigma(-g^{ip}g^{jq}h_{pq}\nabla_if\nabla_jf+
2g^{ij}\nabla_if\nabla_j\frac{\partial f}{\partial s}+\frac{\partial
R}{\partial s})\\
& &\quad \quad \quad \quad \quad \quad \quad \quad  -\
\frac{\partial f}{\partial s}]dV\\
& = & (4\pi\sigma)^{-\frac{n}{2}} \int_M \sigma
e^{-f}g^{ij}h_{ij}(|\nabla f|^2 - \Delta f)dV\\
& & +\ (4\pi\sigma)^{-\frac{n}{2}} \int_M
e^{-f}[\sigma(-g^{ip}g^{jq}h_{pq}\nabla_if\nabla_jf + \frac{\partial
R}{\partial s})-\frac{1}{2}g^{ij}h_{ij}]dV,
\end{eqnarray*}

\hspace{-.6cm}where we note that $$(4\pi\sigma)^{-\frac{n}{2}}
\int_M e^{-f}\cdot 2\sigma g^{ij}\nabla_if\nabla_j\frac{\partial
f}{\partial s}dV = (4\pi\sigma)^{-\frac{n}{2}} \int_M e^{-f}\cdot
2\sigma\frac{\partial f}{\partial s}(|\nabla f|^2-\Delta f)dV,$$
and

\begin{eqnarray*} & & (4\pi\sigma)^{-\frac{n}{2}} \int_M
e^{-f}[\frac{\partial \sigma}{\partial s}(|\nabla f|^2
+R)-\frac{\partial f}{\partial s}]dV\\ & = &
(4\pi\sigma)^{-\frac{n}{2}} \int_M e^{-f}[\frac{\partial
\sigma}{\partial s}(|\nabla f|^2 +R)+\frac{n}{2\sigma}\frac{\partial
\sigma}{\partial s}-\frac{1}{2}g^{ij}h_{ij}]dV\\
& = & (4\pi\sigma)^{-\frac{n}{2}} \int_M
e^{-f}\frac{1}{\sigma}\frac{\partial \sigma}{\partial
s}[\sigma(|\nabla f|^2 +R)+\frac{n}{2}]dV-
(4\pi\sigma)^{-\frac{n}{2}} \int_M
e^{-f}\cdot \frac{1}{2}g^{ij}h_{ij}dV\\
& =& (4\pi\sigma)^{-\frac{n}{2}} \int_M
e^{-f}\frac{1}{\sigma}\frac{\partial \sigma}{\partial
s}[\sigma(2|\nabla f|^2 - 2\Delta
f)+f-\frac{n}{2}+\nu_+]dV-(4\pi\sigma)^{-\frac{n}{2}} \int_M
e^{-f}\cdot \frac{1}{2}g^{ij}h_{ij}dV\\
& = & -(4\pi\sigma)^{-\frac{n}{2}} \int_M e^{-f}\cdot
\frac{1}{2}g^{ij}h_{ij}dV.
\end{eqnarray*}

Since $f(0)=\frac{n}{2}$, we have

\parbox{12cm}
{\begin{eqnarray*} \frac{d^2\nu_+}{ds^2}|_{s=0} & = &
-\frac{1}{Vol(g)}\int_M \sigma
trh\Delta\frac{\partial f}{\partial s}dV\\
&  & +\ \frac{1}{Vol(g)}\int_M(-\frac{n}{2\sigma}\frac{\partial
\sigma}{\partial s}-\frac{\partial f}{\partial
s}+\frac{1}{2}trh)(\sigma \frac{\partial R}{\partial
s}-\frac{1}{2}trh)dV\\
& & +\ \frac{1}{Vol(g)}\int_M(\frac{\partial \sigma}{\partial
s}\frac{\partial R}{\partial s}+\sigma \frac{\partial^2 R}{\partial
s^2}+\frac{1}{2}|h_{ij}|^2)dV.
\end{eqnarray*}}\hfill
\parbox{1cm}{\begin{eqnarray}\end{eqnarray}}

In the following, all quantities are evaluated at $s=0$.\\

Firstly, we have

\begin{eqnarray}
\frac{1}{Vol(g)}\int_M \sigma \frac{\partial^2 R}{\partial s^2}dV &
= & \frac{\sigma}{Vol(g)}\int_M [-\frac{1}{\sigma}|h_{ij}|^2 -
h_{kl}(2\nabla_p\nabla_k
h_{pl}-\Delta h_{kl}- \nabla_k\nabla_l trh) \nonumber\\
& &  \quad \quad \quad \quad \quad -\  \nabla_p [h_{pq}(2\nabla_k h_{kq}- \nabla_q trh)] + \nabla_k (h_{pq}\nabla_k h_{pq})\nonumber\\
& &  \quad \quad \quad \quad \quad +\  \frac{1}{2}\nabla_p trh(2\nabla_k h_{kp} -\nabla_p trh)\nonumber\\
& &  \quad \quad \quad \quad \quad +\   \frac{1}{2}(\nabla_k h_{pq}
\nabla_k h_{pq}- 2\nabla_p h_{kq}\nabla_q h_{kp})]dV\nonumber\\
& = &
\frac{\sigma}{Vol(g)}\int_M[-\frac{1}{\sigma}|h_{ij}|^2-h_{kl}\nabla_p\nabla_k
h_{pl}-\frac{1}{2}|\nabla h|^2-\frac{1}{2}|\nabla
\mathrm{tr}h|^2]dV \nonumber \\
& =&
\frac{\sigma}{Vol(g)}\int_M[-\frac{1}{\sigma}|h_{ij}|^2-h_{kl}(\nabla_k\nabla_p
h_{pl}+R_{kq}h_{ql}+R_{pkql}h_{pq})\nonumber\\
& & \quad \quad \quad \quad \quad -\frac{1}{2}|\nabla
h|^2-\frac{1}{2}|\nabla \mathrm{tr}h|^2]dV\nonumber\\
& = & -\frac{1}{Vol(g)}\int_M \frac{1}{2}|h_{ij}|^2dV \\
& & +\ \frac{\sigma}{Vol(g)}\int_M |\mathrm{div}h|^2 +
Rm(h,h)-\frac{1}{2}|\nabla h|^2 - \frac{1}{2}|\nabla
\mathrm{tr}h|^2dV.\nonumber
\end{eqnarray}

Moreover,
\begin{eqnarray}
\frac{1}{Vol(g)}\int_M \frac{\partial \sigma}{\partial
s}\frac{\partial R}{\partial s}dV& = &
\frac{\sigma}{n}\frac{1}{Vol(g)}\int_M trhdV
\frac{1}{Vol(g)}\int_M \frac{\partial R}{\partial s}dV\nonumber\\
& = & \frac{1}{2n}(\frac{1}{Vol(g)}\int_M trh dV)^2.
\end{eqnarray}

Let $v_h$ be the solution to the following equation,
$$\Delta v_h - \frac{v_h}{2\sigma}=\mathrm{divdiv}h=\nabla_p\nabla_q h_{pq},\ \int_M v_h
=0.$$

Then
\begin{eqnarray*}
& & \frac{1}{Vol(g)}\int_M (-\frac{n}{2\sigma}\frac{\partial
\sigma}{\partial s}-\frac{\partial f}{\partial
s}+\frac{1}{2}trh)\sigma \frac{\partial R}{\partial s}dV\\ & = &
\frac{\sigma}{Vol(g)}\int_M (-\frac{n}{2\sigma}\frac{\partial
\sigma}{\partial s}-\frac{\partial f}{\partial
s}+\frac{1}{2}trh)(\Delta v_h-\frac{v_h}{2\sigma}+
\frac{1}{2\sigma}trh- \Delta trh)dV\\
& = & -\left(\frac{1}{Vol(g)}\int_M \frac{1}{2}trh dV \right)^2 +
\frac{\sigma}{Vol(g)}\int_M v_h(-\Delta \frac{\partial f}{\partial
s}+\frac{1}{2\sigma}\frac{\partial f}{\partial s})dV\\
& &  +\  \frac{\sigma}{Vol(g)}\int_M trh(\Delta \frac{\partial
f}{\partial s}-\frac{1}{2\sigma}\frac{\partial f}{\partial s})dV +
\frac{\sigma}{Vol(g)}\int_M \frac{1}{2}trh(\Delta
v_h-\frac{v_h}{2\sigma}+ \frac{1}{2\sigma}trh- \Delta trh)dV,
\end{eqnarray*}
where we have used (11) to derive the first term in the last equality.\\

Meanwhile,
\begin{eqnarray*}
-\frac{1}{Vol(g)}\int_M \frac{1}{2}trh
(-\frac{n}{2\sigma}\frac{\partial \sigma}{\partial s}-\frac{\partial
f}{\partial s}+\frac{1}{2}trh)& = & -\frac{1}{Vol(g)}\int_M
\frac{1}{2}trh(-2\sigma\Delta\frac{\partial f}{\partial s}-
\sigma\frac{\partial R}{\partial s} + \frac{1}{2}trh).
\end{eqnarray*}
\\
It follows that
\begin{eqnarray*}
& & \frac{1}{Vol(g)}\int_M (-\frac{n}{2\sigma}\frac{\partial
\sigma}{\partial s}-\frac{\partial f}{\partial
s}+\frac{1}{2}trh)(\sigma \frac{\partial R}{\partial
s}-\frac{1}{2}trh)dV\\
& = & \frac{1}{Vol(g)}\int_M \sigma trh\Delta\frac{\partial
f}{\partial s}dV - \frac{1}{Vol(g)}\int_M
\frac{1}{4}(trh)^2dV-\left(\frac{1}{Vol(g)}\int_M \frac{1}{2}trh dV
\right)^2\\
& &  +\  \frac{\sigma}{Vol(g)}\int_M v_h(-\Delta \frac{\partial
f}{\partial s}+\frac{1}{2\sigma}\frac{\partial f}{\partial s})dV+
\frac{\sigma}{Vol(g)}\int_M trh(\Delta \frac{\partial f}{\partial
s}-\frac{1}{2\sigma}\frac{\partial f}{\partial s})dV\\
& & +\ \frac{\sigma}{Vol(g)}\int_M trh(\Delta
v_h-\frac{v_h}{2\sigma}+ \frac{1}{2\sigma}trh- \Delta trh)dV.
\end{eqnarray*}
\\

Now since
\begin{eqnarray*}
\frac{\sigma}{Vol(g)}\int_M v_h(-\Delta \frac{\partial f}{\partial
s}+\frac{1}{2\sigma}\frac{\partial f}{\partial s})dV & = &
\frac{\sigma}{Vol(g)}\int_M v_h(-\frac{n}{4\sigma^2}\frac{\partial
\sigma}{\partial s}+\frac{1}{2}\frac{\partial R}{\partial
s})dV\\
& = & \frac{\sigma}{Vol(g)}\int_M \frac{1}{2}v_h(\Delta v_h
-\frac{v_h}{2\sigma}+ \frac{1}{2\sigma}trh- \Delta trh)dV\\
& = & \frac{\sigma}{Vol(g)}\int_M -\frac{1}{2}|\nabla
v_h|^2-\frac{v_h^2}{4\sigma}+\frac{v_h}{4\sigma}trh-
\frac{1}{2}v_h\Delta trh dV,
\end{eqnarray*}

\hspace{-.6cm}and

\begin{eqnarray*}
\frac{\sigma}{Vol(g)}\int_M trh(\Delta \frac{\partial f}{\partial
s}-\frac{1}{2\sigma}\frac{\partial f}{\partial s})dV& =&
\frac{\sigma}{Vol(g)}\int_M trh(\frac{n}{4\sigma^2}\frac{\partial
\sigma}{\partial s}-\frac{1}{2}\frac{\partial R}{\partial s})dV\\
& = & (\frac{1}{Vol(g)}\int_M\frac{1}{2}trh\ dV)^2\\
&  &  -\ \frac{\sigma}{Vol(g)}\int_M \frac{1}{2}trh(\Delta
v_h-\frac{v_h}{2\sigma}+ \frac{1}{2\sigma}trh - \Delta trh)dV,
\end{eqnarray*}

\hspace{-.6cm}we have

\parbox{10cm}
{\begin{eqnarray*} & & \frac{1}{Vol(g)}\int_M
(-\frac{n}{2\sigma}\frac{\partial \sigma}{\partial s}-\frac{\partial
f}{\partial s}+\frac{1}{2}trh)(\sigma \frac{\partial R}{\partial
s}-\frac{1}{2}trh)dV\\
& = & \frac{1}{Vol(g)}\int_M \sigma trh\Delta\frac{\partial
f}{\partial s}dV+\frac{\sigma}{Vol(g)}\int_M -\frac{1}{2}|\nabla
v_h|^2- \frac{v_h^2}{4\sigma}+ \frac{1}{2}|\nabla trh|^2dV.
\end{eqnarray*}}\hfill
\parbox{1cm}{\begin{eqnarray}\end{eqnarray}}

Substituting (13), (14) and (15) in (12),  we get

\begin{eqnarray*}
\frac{\mathrm{d}^2 \nu_+}{\mathrm{d}s^2}|{s=0} & = & \frac{\sigma}{Vol(g)}\left(\int_M |\mathrm{div}h|^2 + Rm(h,h)-\frac{1}{2}|\nabla h|^2 - \frac{1}{2}|\nabla v_h|^2 - \frac{v_h^2}{4\sigma}dV\right)\\
& &  +\ \frac{1}{2n}\left(\frac{1}{Vol(g)}\int_M trh dV\right)^2\\
&=& \frac{\sigma}{Vol(g)}\int_M <N_+h,h>. \hspace{7cm} \mathrm{Q.E.D.}\\
\end{eqnarray*}

As a simple application, we may briefly discuss the linear stability
of negative Einstein manifolds. Analogue to \cite{Cao}, we say that
a negative Einstein manifold is linearly stable if $N_+\leq 0$,
otherwise it is linearly unstable. As in \cite{Cao}, decompose the
space of symmetric 2-tensors as
$$\mathrm{ker\ div}\oplus \mathrm{im\ div^*},$$
and further decompose $\mathrm{ker\ div}$ as
$$(\mathrm{ker\ div})_0\oplus \mathbb{R}g,$$
where $(\mathrm{ker\ div})_0$ is the space of divergence free
2-tensors $h$ with $\int_M trh=0$. It is easy to see that $N_+$
vanishes on $\mathrm{im\ div^*}$, and on $(\mathrm{ker\ div})_0$
$$N_+=\frac{1}{2}(\Delta_L-\frac{1}{\sigma}),$$
where $\Delta_L=\Delta + 2Rm(\cdot,\cdot)-2Rc$ is the Lichnerowicz
Laplacian on symmetric 2-tensors.\\

Moreover, we may write $(\mathrm{ker\ div})_0$ as
$$(\mathrm{ker\ div})_0=S_0\oplus S_1,$$
where $S_0$ is the subspace of trace free 2-tensors, and $S_1=\{h\in
(\mathrm{ker\ div})_0: h_{ij}=(-\frac{1}{2\sigma}u+\Delta
u)g_{ij}-\nabla_{i}\nabla_{j}u, u\in C^{\infty}(M)\ and \ \int_M
u=0\}$(see e.g. \cite{Bu}).\\

Define $$Tu:=(-\frac{1}{2\sigma}u+\Delta
u)g_{ij}-\nabla_i\nabla_ju.$$ Since $\Delta_L(Tu)=T(\Delta u)$ for
all smooth functions $u$ and $\mathrm{ker}T=\{0\}$, we can see that
the Lichnerowicz Laplacian and the Lapalacian on function space have
the same eigenvalues. Thus $N_+$ is always negative on $S_1$.
Therefore, to study the linear stability of negative Einstein
manifolds, it remains to look at the behavior
of $\Delta_L$ acting on $S_0$ which is the space of transverse traceless 2-tensors.\\

\textbf{Example}\  Suppose that $M$ is an $n$ dimensional compact
real hyperbolic space with $n\geq 3$. By \cite{De} or \cite{Le}, the
biggest eigenvalue of $\Delta_L$ on trace free symmetric 2-tensors
on real hyperbolic space is $-\frac{(n-1)(n-9)}{4}$. Since on $M$ we
have $Rc=-(n-1)g$, $\frac{1}{\sigma}=2(n-1)$. Thus the biggest
eigenvalue of $N_+$ on $S_0$ is not greater than
$-\frac{(n-1)^2}{8}$. It implies that
$M$ is linearly stable for $n\geq 3$.\\

\textbf{Remark}\  When $n=3$, D. Knopf and A. Young (\cite{KnYo})
proved that closed 3-folds with constant negative curvature are
geometrically stable under certain normalized Ricci flow.
R. Ye obtained a more powerful stability result earlier in \cite{Ye}.\\

\textbf{Remark}\  For n=2, R. Hamilton(\cite{Ha}) proved that when
the average scalar curvature is negative, the solution of the
normalized Ricci flow with any initial metric converges to a metric
with constant negative curvature. In particular, they are linearly
stable. On the other hand, in \cite{De2} we see that the biggest
eigenvalue of the Lichnerowicz Laplacian on trace free symmetric
2-tensors is 2. Thus $N_+$ is nonpositive definite on
$(\mathrm{ker\ div})_0$, which also implies the linear stability.\\

\textbf{Remark}  For noncompact case, in \cite{Su}, V. Suneeta
proved certain geometric stability of $\mathbb{H}^n$ using different
methods.

\flushleft
Department of Mathematics, Lehigh University, Bethlehem PA, 18015\\
Email: mez206@lehigh.edu


\begin{thebibliography}{99}
\bibitem{Bu} C. Buzzanca, \textit{The Lichnerowicz Laplacian on
tensors}(Italian), Boll. Un. Mat. Ital. B3 (1984) 531-541

\bibitem{Cao} H.D. Cao, R. Hamilton and T. Ilmanen, \textit{Gaussian densities and stability for some Ricci
solitons}, arXiv:math/0404165v1 [math.DG], 2004

\bibitem{Cao2} H.D. Cao and X.P. Zhu, \textit{A complete proof of the Poincaré and geometrization conjectures---application of the Hamilton-Perelman theory of the Ricci flow}, Asian J. Math. 10 (2006), no. 2, 165-492

\bibitem{De} E. Delay, \textit{Essential spectrum of the lichnerowicz laplacian on two tensor on asymptotically
hyperbolic manifolds}, J. Geom. and Physics 43 (2002), 33-44

\bibitem{De2} E. Delay, \textit{Spectrum of the Lichnerowicz Laplacian on asymptotically hyperbolic
surfaces}, arXiv:0802.3174v1 [math.DG], 2008

\bibitem{Ha} R. Hamilton, \textit{the Ricci flow on surfaces}, Mathematics and general relativity (Santa Cruz, CA, 1986), 237--262
Contemp. Math. 71, Amer. Math. Soc., 1988

\bibitem{KnYo} D. Knopf and A. Young, \textit{Asymptotic Stability of the Cross Curvature flow at a Hyperbolic
Metric}, arXiv:math/0609767v2 [math.DG], 2008

\bibitem{Le} J.M. Lee, \textit{Fredholm Operators and Einstein Metrics on Conformally Compact
Manifolds}, arXiv:math/0105046 [math.DG], 2006

\bibitem{Ni} M. Feldman, T. Ilmanen and Lei Ni, \textit{Entropy and Reduced Distance
for Ricci Expanders}, arXiv:math/0405036v1 [math.DG], 2004

\bibitem{Pe} G. Perelman, \textit{The entropy formula for the Ricci flow and its geometric
applications}, arXiv:math.DG/0211159v1, 2002

\bibitem{Su} V. Suneeta, \textit{Investigating Off-shell Stability of Anti-de Sitter Space in String
Theory}, arXiv:0806.1930v3[hep-th], 2008

\bibitem{Ye} R. Ye, \textit{Ricci Flow, Einstein Metrics and Space
Forms}, Trans. Amer. Math. Soc. 338, no. 2 (1993), 871-896

\end{thebibliography}
\end{document}